\documentclass[a4paper, 11pt]{article}

\usepackage{amsmath}
\usepackage{amsfonts}
\usepackage{amssymb}
\usepackage[english]{babel}
\usepackage{graphicx}
\usepackage{amsthm}
\usepackage{enumerate}
\usepackage{tikz}
\usepackage{bm}

\newtheorem{proposition}{Proposition}

\newtheorem*{cor4.1}{Corollary 4.1}
\newtheorem*{th2.5}{Theorem 2.5}

\theoremstyle{definition}

\date{}

\begin{document}

\title{Some bounds on the number of non isomorphic cyclic $k-$cycle systems of the complete graph }
\author{
Lorenzo Mella \thanks{lorenzo.mella@unibs.it} \ \ \ Gloria Rinaldi  \thanks{Dipartimento di Scienze e Metodi dell'Ingegneria,
Universit\`a di Modena e Reggio Emilia, via Amendola 2, 42122 Reggio
Emilia (Italy) gloria.rinaldi@unimore.it \ \ \ \ \ \ \ \ \ \ \ \ \ \ \  \ \ \ \ Research performed within the activity of INdAM--GNSAGA.}}
%\begin{center}

%\end{center}

\maketitle

\begin{abstract}
\noindent
In a recent paper M. Buratti and M. E. Muzychuck have established some lower bounds on the number of non isomorphic cyclic Steiner Triple Systems systems of order $v\equiv 1$ (mod $6$), \cite{BM}. We complete their results to the case  $v\equiv 3$ (mod $6$). For each odd $k >3$ we also find lower bounds for the number of non isomorphic cyclic  $k$-cycle systems of a complete graph.
\end{abstract}

\noindent \textit{Keywords: cyclic cycle systems; difference systems.}

\noindent\textit{MSC(2010): 05B05; 05B07; 05B10}

\section{Introduction}\label{sec:intro}
Let $K_v$ be the complete graph on $v$ vertices, with $v\ge 3$, and denote by $V$ and $E$ the vertex-set and the edge-set of $K_v$, respectively.
Let $k\ge 3$, a $k$-{\it cycle system} of $K_v$ is a set ${\cal B}$ of $k$-cycles whose edges
partition $E$. ${\cal B}$ is said to be {\it cyclic} if we can identify $V$ with the cyclic additive group $\Bbb Z_v$ and
$B=(b_1, b_2, \dots ,b_k) \in  {\cal B}$ implies that $B+z \in {\cal  B}$ for each $z\in \Bbb Z_v$ and for each $k$-cycle $B$;
where for each $z\in \Bbb Z_v$ the sum $B + z$ means $(b_1+z, b_2+z, \dots , b_k+z)$.

A $k$-cycle system of $K_v$ is also called a $(K_v,C_k)$-{\it design} and in particular a $3$-cycle system of $K_v$ is usually denoted by $STS(v)$ and it is called  a Steiner Triple System of order $v$. The admissible values of $v$ for which there exists a $(K_v,C_k)$-design are those satisfying the following conditions: $v(v-1)\equiv 0$ (mod $2k$) and $v\equiv 1$ (mod $2$).
This is because in a $(K_v,C_k)$-design the number of cycles is  exacly $v(v-1)/2k$, while the number of cycles through any given vertex is $(v-1)/2$.
In what follows, we will assume that these conditions are satisfied. It is easy to prove that these conditions are equivalent to $v\equiv 1$ or $k$ (mod $2k$) whenever $k$ is a prime.

The existence question for $(K_v,C_k)$-designs has been completely settled by Alspach
and Gavlas \cite{AG} in the case of $k$ odd (see also \cite{BD1}) and by Sajna \cite{S} in the even case.

Regarding cyclic $(K_v, C_k)$-designs, existence has been proved whenever $v\equiv  1$
(mod $2k$). This, for $k$ even, was proved in the 1960s by Kotzig \cite{K} and by Rosa \cite{R1,R2}
who also settled the cases of $k = 3, 5, 7$,  \cite{R2} (for the earliest solution of
k = 3 see \cite{P}). The case of k odd and $v\equiv 1$ (mod $2k$) was solved by Buratti and Del Fra \cite{BD1}
and other existence solutions were independently found by Fu and Wu \cite{FW} and by
Bryant et al. \cite{Br}.
The case of $k$ odd and $v\equiv k$ (mod $2k$) was solved by Buratti and Del Fra in \cite{BD1} and  \cite{BD2}, with some exceptions which were solved by Vietri \cite{V}.

To resume, the main results guarantee the existence of a cyclic $(K_v,C_k)$-design whenever $v\equiv 1$ or $k$ (mod $2k$), except for the cases $(v,k)=(9,3)$, $(15,15)$, $(q,q)$, with $q$ a non prime and a prime power. In these cases the non-existence is assured.

%Moreover, the set of values of $v$ for which there exists a cyclic $(K_v,C_k)$-design has
%been completely determined only for $k = 3$, see \cite{P},  and for $k = 4$, see \cite{ZF}.
%In particular a cyclic $(K_v,C_3)$-design, i.e., a cyclic $STS(v)$, exists if and only if $v\equiv 1,6$ (mod $6$)
%and a cyclic $(K_v,C_4)$-design  exists if and only if ?????QUI DEVO CERCARE IL RISULTATO CITATO.

In a cyclic $(K_v,C_k)$-design, given a $k$-cycle $C=(c_1,c_2, \dots , c_k)$ of $K_v$, the {\it type} of $C$ is the cardinality of the stabilizer of $C$ under the action of $\Bbb Z_v$. In particular, the type is $d=1$ if $C+z \ne C$ for each $z\in \Bbb Z_v^* = \Bbb Z_v \setminus \{0\}$, while if $C$ is of type $d>1$ then necessarily $d$ divides both $k$ and $v$ and $C=(c_1, \dots, c_{k/d}, c_1+\ell v/d, \dots, c_{k/d}+\ell v/d, \dots, c_1+\ell (d-1)v/d, \dots, c_{k/d}+\ell (d-1)v/d)$, with $\ell$ and $d$ coprime. In both these cases, the list of partial differences from C is the multiset $\partial C = \{\pm (c_{i+1}-c_i): 1\le i \le k/d -1\}\cup \{\pm (c_1-c_{k/d}+\ell v/d)\}$.

More generally, if ${\cal F}= \{C_1, C_2, \dots , C_n\}$ is a set of $k$-cycles, the list of partial
differences from ${\cal F}$ is the multiset $\partial {\cal F}= \bigcup _i \partial C_i$.

Note that if $C$ is a cycle of type 1, then $\partial C$ contains $2k$ distinct elements and it is the list $\Delta C$ of differences from $C$ in the usual sense.

The orbit of a cycle $C$ will be denoted by $Orb_{\Bbb Z_v}(C)$. Obviously, if $C$ is of type $d$ its length orbit is $v/d$.

As a consequence of \cite{B2} it is well known that if  ${\cal F}= \{C_1, C_2, \dots , C_n\}$ is a set of $k$-cycles, and $d_i$ is the type of $C_i$, $1\le i\le n$,  then the cycles of the set $\bigcup_{i=1}^n Orb_{\Bbb Z_v}(C_i)$ form a cyclic $(K_v,C_k)$-design if and only if $\partial {\cal F}$ covers $\Bbb Z_v^*$ exactly once. If this is the case, ${\cal F}$ is called a $(K_v,C_k)$-{\it difference system}, or simply a {\it difference system} when $k$ and $v$ are clear from the context.  Each cycle $C_i$ of ${\cal F}$ is called a {\it starter cycle}.
Moreover, we can always suppose that each starter cycle contains $0$ as a vertex.

If $v \equiv 1$ (mod $2k$) then necessarily each starter cycle is of type 1, while if $v\equiv k$ (mod $2k$) a difference system can contain starter cycles of type $d_i > 1$ with $d_i$ dividing both $v$ and $k$.

In this paper we will also consider the list $\Lambda ({\cal F})$ of all the lists of partial differences of the cycles of ${\cal F}$:

$$\Lambda ({\cal F})=\{\partial C_1, \dots, \partial C_n\}.$$

%Recall that a $2-(v,k,\lambda)$ design is a a pair $(V,{\cal B})$ where $V$ is a set of $v$ {\it points} and ${\cal B}$ is a set of $k$-subsets ({\it blocks}) of $V$ such that any $2$-subset of $V$ is contained in exactly $\lambda$ blocks (see \cite{BJL} and \cite{Hand} for a general background). Such a design is said to be cyclic if $V=\Bbb Z_v$ and ${\cal B}$ is invariant under the addition of $\Bbb Z_v$.
%Two designs $(V,{\cal B})$ and $( V',{\cal B'})$ are said to be isomorphic if there is a bijection $\alpha: V\longrightarrow V'$ turning ${\cal B}$ into ${\cal B'}$. To establish the number of $2-(v,k,\lambda)$ designs up to isomorphism is in general not feasible and very little is known about the number of cyclic $2-(v,k,\lambda)$-designs up to isomorphism, indeed it is also hard to establish if a single $2-(v,k,\lambda)$ design exists.

Two $(K_v,C_k)$-designs, with $k$-cycle sets ${\cal B}$ and ${\cal B'}$ respectively, are said to be {\it isomorphic} if there exists a permutation of $V$ turning ${\cal B}$ into ${\cal B'}$. In what follows we will denote by $NC(v,k)$ the number of cyclic $(K_v,C_k)$-designs up to isomorphism. We will denote by $\phi$ the Euler's totient function.

As a consequence of the famous Bays-Lambossy Theorem, see \cite{Ba} and its generalizations proved by P\'{a}lify, \cite{Pa}, and Muzychuck, \cite{M}, the following result can be stated. We omit its proof since it is a repetition of that of Corollary 1.3 of \cite{BM}.

%In \cite{BM} it is proved that the existence of just one cyclic $(K_v,C_k)$-design necessarily implies
%$NC(v,k) > \lceil \frac{2^n}{\phi(v)}\rceil$ with $n=\frac{v-1}{k^2-k}$ or $n=\frac{v-k}{k^2-k}$ according to the congruence class of $v$ modulo $k^2-k$. In the same paper it is proved that that $NC(v,3)\ge 2$ for all admissible values of $v\ge 481$, and $NC(6n+1,3)>2.49^n$ for any $n\not\equiv 4$ (mod $7$), and $NC(6n+1,3)>3.35^n$ for any $n\equiv 4$ (mod $7$) sufficiently large. These results are better than those given in \cite{CR}. To obtain their results, Buratti and Muzychuk used the following Lemma \ref{L1}, which is a consequence of the famous Bays-Lambossy Theorem, see \cite{Ba} and its generalizations proved by P\'{a}lify, \cite{Pa}, and Muzychuck, \cite{M}.

\begin{proposition}\label{L1}
If {\cal D} is a set of pairwise distinct cyclic $(K_v,C_k)$-designs, then the relatio  $NC(v,k)\ge \lceil\frac{|{\cal D}|}{\phi(v)}\rceil$ holds.
\end{proposition}

If we want to establish some lower bounds for $NC(v,k)$, it  could be important to produce ``many" distinct cyclic $(K_v,C_k)$-designs.
With this aim, we briefly explain the proof of the following Proposition \ref{L2} using a procedure similar to the one used in  Theorem 2.1 of \cite{BM}.
First of all, for each $k$-cycle $C=(c_1, c_2, \dots, c_k)$ we denote by $-C$ the $k$-cycle $(-c_1,-c_2, \dots, -c_k)$.

\begin{proposition}\label{L2}
Let ${\cal F}= \{C_1, C_2, \dots , C_n\}$ be a $(K_v,C_k)$-difference system and let $1\le s \le n$ be the number of cycles of ${\cal F}$ which are of type 1. We have $NC(v,k)\ge \lceil\frac{2^s}{\phi(v)}\rceil$.
\end{proposition}

\noindent
{\em Proof.} Without loss in generality, suppose $C_1, \dots, C_s$ to be the cycles of ${\cal F}$ which are of type 1. Observe that $C_i$ and $-C_{i}$ give the same partial difference set, i.e. $\partial (-C_i) = \partial C_i$. Moreover, all the difference of $C_i$ are distincts, i.e. $\partial C_i = \Delta C_i$.
Take the $s$-tuples $\sigma=(\sigma_1, \dots, \sigma_s) \in \{1,-1\}^s$ and let ${\cal F}_{\sigma}=\{\sigma_1C_1, \dots , \sigma_s C_s\}\cup \{C_{s+1}, \dots, C_n\}$. It is obvious that ${\cal F}_{\sigma}$ is a difference system itself. Let $D_{\sigma}$ be the cyclic $(K_v,C_k)$-design obtained by the orbits of the starter cycles of ${\cal F}_{\sigma}$. We show that ${\cal D}=\{D_{\sigma} \ | \ \sigma \in \{1,-1\}^s\}$ has no repeated $(K_v,C_k)$-designs. In fact, assume that $\sigma$ and $\sigma'$ are distinct $s$-tuples of $\{1,-1\}^s$ and, without loss of generality assume that $\sigma_i=1$ and $\sigma'_i = -1$ for a suitable $i$, $1\le i \le s$. Observe that $C_i$ is neither a translate of $\sigma'_jC_j$ for $i\ne j$, $1\le j \le s$, nor a translate of $C_j$, $s+1 \le j \le n$,  otherwise it should be $\partial C_i = \partial C_j$ contradicting the fact that $\partial {\cal F}$ has no repeated elements.
Now assume that $C_i$ is a translate of $\sigma'_iC_i = -C_i$, say $C_i= -C_i +t$. Since $0$ is a vertex of $C_i=(0,c_2, \dots, c_k)$, and $[0,c_2]$ is a unique edge with difference set $\{\pm c_2\}$, we necessarily have $[0,c_2]=[t,-c_2+t]$ and $[0,c_k]=[t,-c_k+t]$. This implies either $t=0$ or $t=c_2$. In the first case we should have $c_2=-c_2$: a contradiction as $v$ is odd; in the second case we should have $c_2=c_k$ and again a contradiction.

We conclude that the cycle $C_i$ of
$D_{\sigma}$ cannot be in $D_{\sigma'}$, hence $D_{\sigma}$ and $D_{\sigma'}$ are different. The assertion then follows from Proposition \ref{L1}. $\qed$

\vskip 0.5truecm
\noindent
We will apply these results in the next sections.

We recall that when $k=3$ and $v\equiv 1$ (mod $6$) a lower bound of $NC(6n+1,3)$ was obtained in \cite{BM}, while the case $v\equiv 3$ (mod $6$) was left open. We complete the results of \cite{BM} and determine lower bounds of $NC(6n+6,3)$. We also consider the case $k=5$ a part. We determine lower bounds of $NC(v,5)$, whenever both the necessary conditions $v\equiv 1,5$ (mod $10$) hold. In a similar manner, we provide some lower bounds of $NC(v,k)$, for $k$ odd and $k > 5$.

In the rest of the paper we will make use of Skolem-type sequences, whose definitions are given below.

A {\it Skolem sequence} of order $n$ is a sequence of $n$ integers $(s_1, \dots,s_n)$ such that
$\bigcup_{i=1}^{n}\{s_i,s_i+i\}=\{1,2,\dots, 2n+1\} - \{k\}$ where either $k=2n+1$ or $k=2n$ according to whether $n\equiv 0,1$ (mod $4$), or $n\equiv 2,3$ (mod $4$). In the last case one usually speaks of a {\it hooked Skolem sequence}.
For the rest of the paper, we will simply talk about a Skolem sequence according to the appropriate residue class of $n$. It is well known that a Skolem sequence of order $n$ exists for every positive $n$, \cite{Hand}. Moreover, it was known for a long time, \cite{A}, that $Sk(n)$, the number of distinct Skolem sequences of order $n$, is bounded from below by $2^{\lfloor n/3 \rfloor}$. A tighter lower bound was obtained in \cite{DG} at Corollaries 3.3. and 3.4. It states that $Sk(n)\ge (6.492)^{(2n-1)/7}$ for sufficiently large $n \equiv 4$ (mod $7$).

A {\it split Skolem sequence} of order $n$ is a sequence of $n$ integers $(s_1,s_2,\dots , s_n)$ such that $\bigcup_{i=1}^{n}\{s_i,s_i+i\}=\{1,2,\dots, 2n+2\} - \{n+1, k\}$ with either $k=2n+2$ or $k=2n+1$ according to whether $n\equiv 0,3$ (mod $4$) or $n\equiv 1,2$ (mod $4$). In the last case one usually speaks of a {\it split-hooked Skolem sequence}. Split Skolem sequences and split-hooked Skolem sequences are also known as {\it Rosa sequences} and {\it hooked Rosa sequences}, respectively. In what follows we will simply talk about a split Skolem sequence (or a Rosa sequence) according to the appropriate residue class of $n$. It is well known that a split Skolem sequence of order $n$ exists for every positive $n>1$,  \cite{Hand}.

Again use the notation $Sk(n)$ for the number of distinct split Skolem sequences of order $n$. It was known for a long time, \cite{A}, \cite{BGG}, that $Sk(n)$  is essentially bounded below by $2^{\lfloor n/3 \rfloor}$. A more tight lower bound was obtained in \cite{DG} at Corollaries 3.2 and 3.5. They states that for sufficiently large $t\equiv 0,3$ (mod $4$), there are more than $(6.492)^{2t+1}$ split Skolem sequences of order $7t+3$ and more than $(6.492)^{2t+1}$ split-hooked Skolem sequences of order $7t+5$.

\section{Lower bounds of $NC(6n+3,3)$}

Recall that a $(K_v, C_3)$-design, i.e. a $STS(v)$,  exists if  $v\equiv 1, 3$ (mod $6$). In other words, if it is either $v=6n+1$ or $v=6n+ 3$. In both cases, the existence of a cyclic $STS(v)$ is assured for each $n\ge 1$ but $v=9$.
Some lower bounds for $NC(6n+1,3)$ were established in Proposition 3.1 of \cite{BM}. In this section we complete the result and give some lower bounds for $NC(6n+3,3)$.

Let $S=(s_1,\dots, s_n)$ be a split Skolem sequence of order $n$, or a hooked-split Skolem sequence depending on the appropriate residue class of $n$. The set ${\cal F}= \{C_1, C_2, \dots , C_n\}\cup \{B\}$,  where $C_i =(0,i,s_i+i+n)$, $i=1, \dots, n$ and $B=(0,2n+1,4n+2)$ is a difference system which determines a cyclic $STS(6n+3)$. In fact we have:
$$\partial{\cal F}=\{\pm i, \pm(s_i+n), \pm(s_i+i+n), i=1,\dots, n\}\cup \{\pm(2n+1)\}= \Bbb Z_{6n+3}^*.$$
\noindent
Moreover,  each starter cycle $C_i$ is of type $1$, while the starter cycle $B$ is of type $3$.

We can prove the following 

\begin{proposition}\label{Gen}
Let $NC(6n+3,3)$ be the number of cyclic non-isomorphic $STS(6n+3)$s. Therefore:

\begin{itemize}

\item If $n\equiv 3, 5, 24, 26 \  (mod \ 28)$ is sufficiently large,
\par\noindent
then $NC(6n+3,3)\ge (3.35)^n$.

\item If $n\not\equiv 3, 5, 24, 26 \  (mod \ 28)$ is sufficiently large,
\par\noindent
then $NC(6n+3,3)\ge (2.49)^n$.

\end{itemize}

\end{proposition}

\noindent
{\em Proof.} Let ${\cal F }$ be the difference system above described. Take the $n$-tuples $\sigma=(\sigma_1, \dots, \sigma_n) \in \{1,-1\}^n$ and let ${\cal F}_{\sigma}=\{\sigma_1C_1, \dots , \sigma_n C_n\}\cup \{B\}$. Let $D_{\sigma}$ be the cyclic $STS(6n+3)$ obtained by the orbits of the starter cycles of ${\cal F}_{\sigma}$. We know from Proposition \ref{L2} that ${\cal D}=\{D_{\sigma} \ | \ \sigma \in \{1,-1\}^n\}$ contains $2^n$ different $STS(6n+3)$s. Now observe that if $S$ and $S'$ are two different split Skolem sequences, respectively split-hooked Skolem sequences, of the same order $n$ and if ${\cal F}$ and ${\cal F}$' are the related difference systems giving rise to the cyclic $STS(6n+3)$s $D$ and $D'$ respectively, then $D\ne D'$ since $\Lambda ({\cal F}) \ne \Lambda ({\cal F}')$. We can conclude that the number of distinct cyclic $STS(6n+3)$s is at least $2^n \times Sk(n)$, where $Sk(n)$ is the number of distinct split Skolem sequences, respectively split-hooked Skolem sequences, of order $n$ according to whether $n\equiv 0,3$ (mod $4$) or $n\equiv 1,2$ (mod $4$), $n\ne 1$. If $n\equiv 0,3$ (mod $4$) and $n=7t+3$, then  for sufficiently large $t\equiv 0,3$ (mod $4$) we have: $n\equiv 3,24$ (mod $28$), $2t+1= (2n+1)/7$ and $Sk(n) > (6.492)^{(2n+1)/7}$.
If $n\equiv 1,2$ (mod $4$), $n\ne 1$ and $n=7t+5$, then for sufficiently large $t\equiv 0,3$ (mod $4$) we have: $n\equiv 5,26$ (mod $28$), $2t+1= (2n-3)/7$ and $Sk(n) > (6.492)^{(2n-3)/7}$.
We can state:

\noindent
If $n\equiv 3,24$ (mod $28$), then:
$$NC(6n+3) > \frac{2^n \times (6.492)^{(2n+1)/7}}{\phi(6n+3)} \ge \frac{2^n \times (6.492)^{(2n+1)/7}}{6n+2}$$

\noindent
If $n\equiv 5,26$ (mod $28$), then:

$$NC(6n+3) > \frac{2^n \times (6.492)^{(2n-3)/7}}{\phi(6n+3)} \ge \frac{2^n \times (6.492)^{(2n-3)/7}}{6n+2}$$

Numeric calculations show that the fractions above are both bigger than $(3.35)^n$ starting from $n\ge 444$.

In all the other cases we use the bound  $\frac{2^n \times 2^{\lfloor n/3 \rfloor}}{\phi(6n+3)} \ge \frac{2^n \times 2^{\lfloor n/3 \rfloor}}{6n+2}$. In this case numeric calculations show that this fraction is bigger than $(2.49)^n$ starting from $n\ge 702$.
$\qed$

\section{Lower bounds of $NC(10n+k,5)$, $k=1,5$}
In what follows we generalize the previous Proposition \ref{L2}.

\begin{proposition}\label{L3}
Let ${\cal F}=\{C_1,\dots C_n\}$ be a $(K_v,C_5)$-difference system and let $1\le s\le n$ be the number of cycles of ${\cal F}$ which are of type 1.
We have $NC(v,5)\ge \lceil\frac{24^s}{\phi(v)}\rceil$.
\end{proposition}
\noindent
{\em Proof.} Let $C$ be a type 1 $5$-cycle of ${\cal F}$ with $C=(0,a_1,a_2,a_3,a_4)$. Set $d_i=a_i - a_{i-1}$, with $i=1,\dots, 5$ and $a_5=0$.
We can also write  

$$C=(0, d_1, \sum_{i=1}^{2}d_i,\sum_{i=1}^{3}d_i,\sum_{i=1}^{4}d_i).$$

\noindent
Let $\overrightarrow{D_C}=(d_1,d_2,d_3,d_4,d_5)$ be the oriented $5$-cycle associated to $C$. If we consider a permutation of the vertices of $\overrightarrow{D_C}$ we obtain an oriented $5$-cycle with an associate type 1 cycle of $K_v$. 
Namely, if $\sigma$ is a permutation of $S_5$, the symmetric group on $5$ elements,   and we denote by $\sigma_i$  the image of $i=1,\dots, 5$ under  $\sigma$, we construct:

$$\overrightarrow{D_{C_{\sigma}}}=(d_{\sigma_1},d_{\sigma_2},d_{\sigma_3},d_{\sigma_4},d_{\sigma_5}) \ \ \ C_{\sigma}=(0, d_{\sigma_1}, \sum_{i=1}^{2}d_{\sigma_i},\sum_{i=1}^{3}d_{\sigma_i},\sum_{i=1}^{4}d_{\sigma_i}).$$ 

It is easy to verify that $C_{\sigma}$ is still a $5$-cycle and it is of type 1 since $\partial C = \partial C_{\sigma}$.

Let $\tau \in S_5$, construct: 

$$\overrightarrow{D_{C_{\tau}}}=(d_{\tau_1},d_{\tau_2},d_{\tau_3},d_{\tau_4},d_{\tau_5}) \ \ \ 
C_{\tau}=(0, d_{\tau_1},\sum_{i=1}^{2} d_{\tau_i},\sum_{i=1}^{3} d_{\tau_i},\sum_{i=1}^{4} d_{\tau_i}).$$

Observe that $\partial C= \partial C_{\sigma}=\partial C_{\tau}$. Furthermore $C_{\sigma}\in Orb_{Z_v}(C_{\tau})$ if and only if $\overrightarrow{D_{C_{\sigma}}}=\overrightarrow{D_{C_{\tau}}}$. In fact, if $C_{\sigma}$ and $C_{\tau}$ are in the same orbit, then their associated oriented sequences of the elements $\{d_1, \dots, d_5\}$ must be the same. Viceversa, if $\overrightarrow{D_{C_{\sigma}}}=\overrightarrow{D_{C_{\tau}}}$ with $d_{\tau_r}=d_{\sigma_1}$, it is $C_{\sigma}+t=C_{\tau}$, with either $t=0$ or $t=\sum_{j=2}^{r}d_{\sigma_{5-r+j}}$ according to whether $r=1$ or $2\le r \le 5$. We conclude that for each type 1 $5$-cycle $C$ of ${\cal F}$, the set $\{C_{\alpha} \ | \ \alpha \in S_5 \}$ can be partitioned into $4!$ classes, where each class contains $5$ cycles in the same orbit, while  two cycles not in the same class have different orbits.
We have  ${\cal F}=\{C_1, \dots , C_n\}$ and without loss in generality, we can suppose that  $C_1, \dots, C_s$ are of type 1.

Denote by ${\cal A}_i=\{{\cal A}_i^1, {\cal A}_i^2,\dots, {\cal A}_i^{4!}\}$ the set formed by the $4!$ distinct classes of cycles arising from $C_i$, $1\le i \le s$, as above described. 

Let $B=\{1,2,\dots, 4!\}^s$. For each $\varphi=(r_1, \dots, r_s)\in B$ 
\noindent
let ${\cal F}_{\varphi}=\{C_1^{r_1}, \dots, C_s^{r_s},C_{s+1}, \dots, C_n\}$ with $C_i^{r_i}\in {\cal A}_i^{r_i}$. Obviously ${\cal F}_{\varphi}$ is a difference system. Denote by ${\cal D}_{\varphi}$ the cyclic $(K_v,C_5)$-design with difference set ${\cal F}_{\varphi}$. If $\psi \in B$ and $\psi \ne \varphi$, then ${\cal D}_{\varphi}\ne {\cal D}_{\psi}$ since they contain cycles with different orbits. Since $B$ contains $(4!)^s$ elements, we have constructed $(4!)^s$ distinct designs and the assertion follows 
from Proposition \ref{L1}. $\qed$

\vskip0.3truecm\noindent
{\bf Remark.}
Let $C$ be the $5$-cycle $C=(0,a_1,a_2,a_3,a_4)$ with associated oriented $5$-cycle $\overrightarrow{D_C}=(d_1,d_2,d_3,d_4,d_5)$. We can observe that the $5$-cycle $-C$ is obtained from $\overrightarrow{D_{C_{\alpha}}}=(d_5,d_4,d_3,d_2,d_1)$, with $\alpha \in S_5$, $\alpha =(15)(24)(3)$. Therefore, Proposition \ref{L3} is a generalization of Proposition \ref{L2}. In fact, if we apply Proposition \ref{L2} to a $(K_v,C_5)$-difference system ${\cal F}$, we just consider $2$ of the $4!$ classes arising from each cycle $C$ of ${\cal F}$.

Proposition \ref{L3} fails if we try to prove it for cyclic $(K_v, C_k)$-designs with $k>5$. In fact, if $C$ is a type 1 $k$-cycle with $\overrightarrow{D_C}=(d_1,d_2, \dots, d_k)$ and $\alpha \in S_k$, then $C_{\alpha}=(0,d_{\alpha_1}, d_{\alpha_1}+d_{\alpha_2}, \dots, \sum_{i=1}^{k-1}d_{\alpha_i})$ is not necessarily a $k$-cycle.

Finally, if we consider cyclic $(K_v, C_3)$-designs, then for each type 1 cycle $C=(0, a_1,a_2)$ of a difference system ${\cal F}$, we have exacly $2$ classes arising from $C$: the one determined by $\overrightarrow{D_C}=(d_1,d_2,d_3)$ and that determined by $\overrightarrow{D_{C_{\alpha}}}=(d_3,d_2,d_1)$, with $\alpha \in S_3$, $\alpha=(13)(2)$ and $C_{\alpha}=-C$. Therefore Proposition \ref{L2} and Proposition \ref{L3} are exactly the same in this case.

\subsection{The case $NC(10n+1,5)$}

The construction of a starter family ${\cal F}$ of a cyclic $(K_{10n+1}, C_5)$-design is given in \cite{BD1}. More precisely, if $n=1$ we have ${\cal F}=\{C\}$ with $C=(0,-3,-4,3,-6)$. If $n>1$,  let $S=(s_1, \dots,s_n)$ be a Skolem sequence  of order $n$, then we have ${\cal F}=\{C_1, \dots, C_n\}$, where $C_i =(0,s_i+i,i,-2n,i+3n)$ with $C_1$ replaced by $C'_1=(0,s_1+1,1,5n+1,2n)$ if the Skolem sequence is a hooked Skolem sequence, i.e. in the case $n\equiv 2,3$ (mod $4$).

The known bounds on the number $Sk(n)$ of Skolem sequences yield the following

\begin{proposition}\label{Gen5}
Let $n >1$ and let $NC(10n+1,5)$ be the number of cyclic non-isomorphic $(K_{10n+1},C_5)$-designs. Therefore:

\begin{itemize}

\item If $n\not\equiv 4$ (mod $7$) is sufficiently large, then
$NC(10n+1,5)\ge (29.7)^n$.

\item If $n \equiv 4$ (mod $7$) is sufficiently large, then
$NC(10n+1,5)\ge (38.4)^n$.

\end{itemize}
\end{proposition}

\noindent
{\em Proof.} Let ${\cal F }$ be the $(K_v,C_5)$-difference system obtained in \cite {BD1} and above described. It is obvious that all the cycles of ${\cal F}$ have trivial stabilizer in $\Bbb Z_{10n+1}$, in fact  $1$ is the unique common divisor of $5$ and $10n+1$. Let ${\cal D}$ be the cyclic design arising from ${\cal F}$.  Following the notation of Proposition \ref{L3}, let $D_{\varphi}$ be the cyclic $(K_{10n+1}, C_5)$-design obtained by the orbits of the starter cycles of ${\cal F}_{\varphi}$. We know from Proposition \ref{L3} that $\varphi$ can be chosen in $(4!)^n$ manners so that  a set of $\lceil\frac{(4!)^n}{\phi(10n+1)}\rceil$  different designs arise from $\cal F$.  Now observe that if $S$ and $S'$ are two different Skolem sequences of the same order $n$ and if ${\cal F}$ and ${\cal F}$' are the related difference systems giving rise to the cyclic $(K_v,C_5)$-designs $D$ and $D'$ respectively, then $D\ne D'$ since $\Lambda ({\cal F}) \ne \Lambda ({\cal F}')$. We can conclude that the number of distinct cyclic $(K_v,C_5)$-designs is at least $\dfrac{(4!)^n\cdot Sk(n)}{\phi(10n+1)}$ where $Sk(n)$ is the number of distinct Skolem sequences, respectively hooked Skolem sequences, of order $n$ according to whether $n\equiv 0,1$ (mod $4$) or $n\equiv 2,3$ (mod $4$). 
We known that $Sk(n)\ge 2^{\lfloor n/3 \rfloor}$ whenever $n\not \equiv 4$ (mod $7$), while $Sk(n)\ge (6.492)^{(2n-1)/7}$ whenever $n\equiv 4$ (mod $7$). 
In the first case, we obtain:
$NC(10n+1,5)\ge \dfrac{(24)^n\cdot 2^{\lfloor n/3 \rfloor -1}}{5n}$ and an easy calculation shows that this number is bigger than $(24)^n$ for $n \ge 27$ and it is bigger than $(29.76)^n$ for $n \ge 570$.

In the second case, we obtain $NC(10n+1,5)\ge \dfrac{(24)^n \cdot (6.492)^{(2n-1)/7}}{10n}$. An easy calculation shows that this number is bigger than $(31.35)^n$ for $n\ge 25$ and it is bigger than $(38.4)^n$ for $n \ge 116$.  $\qed$

\vskip0.5truecm
If we consider the case $n=1$, the previous Proposition \ref{L3} simply says that $NC(11,5) \ge 3$.

\subsection{The case $NC(10n+5,5)$}
The construction of a starter system ${\cal F}$ of a cyclic $(K_{10n+5}, C_5)$-design is given in \cite{BD1}. More precisely, set $m=2n+1$ and let $H=m\Bbb Z_{5m}=\{0,m,2m,3m,4m\}$ be the subgroup of the additive group $\Bbb Z_{5m}$. Consider the multipartite graph $K_{m\times 5}$ whose parts are the cosets $H+h$, $0\le h\le m-1$. A difference system ${\cal F'}$ whose orbit under $\Bbb Z_{5m}$ gives a cyclic $(K_{m\times 5},C_5)$-design is constructed as follows, see Lemma 3.1 and Theorem 3.2 of \cite{BD1}. If $n=1$, then ${\cal F'}= \{B\}$, with $B=(0,-1,1,-6,4)$. While, if $n>1$, then ${\cal F'}= \{B_1, \dots, B_n\}$, with $B_i=(b_{i1},b_{i2},b_{i3},b_{i4,},b_{i5})$ defined by the following rules:

$$b_{ij}: \left\{\begin{array} {c} \dfrac{m(j-1)}{2} \ \ \  \text{for} \  j \ \text{odd,} \ \ j\ne 5  \\ m\left(2 - \dfrac{j}{2}\right)-i \ \ \  \text{for} \  j \ \text{even,} \end{array}  \right.$$

$$b_{i5}=r_i +(2m-n-1),$$

\noindent
where $(r_1,\dots, r_n)$ is a either a split Skolem sequence or a hooked split Skolem sequence of order $n$ according to whether $n\equiv 0,3$ or $n\equiv 1,2$ (mod $4$). The set ${\cal F'}$ is a difference system for a cyclic $(K_{m\times 5},C_5)$ design as $\partial {\cal F'} = \Bbb Z_{5m} - m\Bbb Z_{5m}$, Theorem 3.2 of \cite{BD1}.

Now consider the cyclic $(K_5,C_5)$-design with vertex set $\Bbb Z_5$ and cycle set given by $\{\overline{B}_1,\overline{B}_2\}$, with $\overline{B}_1=(0,1,2,3,4)$ and $\overline{B}_2=(0,2,4,1,3)$, see Theorem 4.2 of \cite{BD1}.
Set $m\overline{B}_1=(0,m,2m,3m,4m)$ and $m\overline{B}_2=(0,2m,4m,m,3m)$. Applying Theorem 4.1 of \cite{BD1}, we can state that the set
$\{m\overline{B}_1,m\overline{B}_2\}\cup {\cal F'}$ is a $(K_{10n+5},C_5)$ difference system.  
Moreover,  $m\overline{B}_1$ and $m\overline{B}_2$ are of type $5$, with $\partial m\overline{B}_1 = \pm\{m\}$ and $\partial m\overline{B}_2 = \pm\{2m\}$, while all the cycles of ${\cal F'}$ are of type 1. In fact, since $\partial {\cal F'} = \Bbb Z_{5m} - m\Bbb Z_{5m}$, each cycle of ${\cal F'}$ gives rise to $10$ distinct differences of $\Bbb Z_{5m}$.

\begin{proposition}\label{Gen5.5}
Let $n > 1$ and let $NC(10n+5,5)$ be the number of cyclic non-isomorphic $(K_{10n+5},C_5)$-design. Therefore:

\begin{itemize}

\item If $n\equiv 3, 5, 24, 26 \  (mod \ 28)$ is sufficiently large,
\par\noindent
then $NC(10n+5,5)\ge (38.4)^n$.

\item If $n \not \equiv 3, 5, 24, 26 \  (mod \ 28)$ is sufficiently large,
\par\noindent
then $NC(10n+5,5)\ge (29.7)^n$.

\end{itemize}

\end{proposition}

\noindent
{\em Proof.} Let ${\cal F }$ be the $(K_{10n+5},C_5)$-difference system above described. Let ${\cal D}$ be the cyclic design arising from ${\cal F}$.  Following the notation of Proposition \ref{L3}, let $D_{\varphi}$ be the cyclic $(K_{10n+5}, C_5)$-design obtained by the orbits of the starter cycles of ${\cal F}_{\varphi}$. Since $n$ is the number of cycles of ${\cal F}$ with trivial stabilizer, we know from Proposition \ref{L3} that $\varphi$ can be chosen in $(4!)^n$ manners so that  a set of $\lceil\frac{(4!)^n}{\phi(10n+5)}\rceil$  different designs arise from $\cal F$.  Now observe that if $n >1$,  $S$ and $S'$ are two different split or hooked split sequences of the same order $n$, and if ${\cal F}$ and ${\cal F}$' are the related difference systems giving rise to the cyclic $(K_{10n+5},C_5)$-designs $D$ and $D'$ respectively, then $D\ne D'$ since $\Lambda ({\cal F}) \ne \Lambda ({\cal F}')$. We can conclude that the number of distinct cyclic $(K_{10n+5},C_5)$-designs is at least $\dfrac{(4!)^n\cdot Sk(n)}{\phi(10n+1)}$ where $Sk(n)$ is the number of distinct split sequences, respectively hooked split sequences, of order $n$ according to whether $n\equiv 0,1$ (mod $4$) or $n\equiv 2,3$ (mod $4$). We can use the known bounds for $Sk(n)$ given in \cite{DG} to obtain the result.
Namely, if $n\equiv 0,3$ (mod $4$) and $n=7t+3$, then  for sufficiently large $t\equiv 0,3$ (mod $4$) we have: $n\equiv 3,24$ (mod $28$), $2t+1= (2n+1)/7$ and $Sk(n) > (6.492)^{(2n+1)/7}$.
If $n\equiv 1,2$ (mod $4$), $n\ne 1$ and $n=7t+5$, then for sufficiently large $t\equiv 0,3$ (mod $4$) we have: $n\equiv 5,26$ (mod $28$), $2t+1= (2n-3)/7$ and $Sk(n) > (6.492)^{(2n-3)/7}$. Therefore, if $n$ is sufficiently large and either $n \equiv 3, 24$ (mod $28$) or $n \equiv 5, 26$ (mod $28$), we have respectively $NC(10n+5,5)\ge \dfrac{(24)^n \cdot (6.492)^{\frac{2n+1}{7}}}{\phi(10n+5)}$ or $NC(10n+5,5)\ge \dfrac{(24)^n \cdot (6.492)^{\frac{2n-3}{7}}}{\phi(10n+5)}$. An easy calculation shows that both these numbers are bigger than $(31.35)^n$ if $n \ge 25$ and are bigger than $(38.4)^n$ if $n\ge 136$. If $n \not \equiv 3, 5, 24, 26$ (mod $28$) is sufficiently large, then $NC(10n+5,5)\ge \dfrac{(24)^n \cdot 2^{\lfloor n/3\rfloor}}{\phi(10n+5)}$ and this number is bigger than  $(24)^n$ when $n\ge 27$ and it is bigger than $(29.76)^n$ if $n \ge 570$. $\qed$

\vskip0.5truecm
If we consider the case $n=1$, the previous Proposition \ref{Gen5.5} simply says that $NC(15,5) \ge 3$.

\section{Lower bounds of $NC(2nk+1,k)$ and $NC(2nk+k,k)$ with $k>5$ odd}

In this section we obtain lower bounds for the number of non isomorphic cyclic  $(K_v,C_k)$-designs when $k>5$, $k$ odd (except for $v=3k$ and $v=2k+1$). For this purpose we use the examples of cyclic $(K_v,C_k)$-designs, with either $v=2nk+1$ or $v=2nk+k$, $n >1$, which are obtained in \cite{BD1} \cite{BD2}, \cite{V} via Skolem-type sequences. We will use the same technique already applied above and in \cite{BM} and the obtained lower bounds  depends upon both $n$ and $k$. They are significative, in the sense of high,  if $k$ itself is not too high, in particular if $k$ is linear with $n$.

First of all we recall the contructions of \cite{BD1}.  More precisely, suppose $n\equiv 0,1$ (mod $4$) and fix a Skolem sequence $(s_1, \dots, s_n)$ of order $n$. For $i=1, \dots , n$, consider the $k$-cycle defined as follows:
\vskip0.3truecm\noindent
$b_{i1}=0$, \ \ $b_{i2}=-s_i$

$$b_{ij}: \left\{\begin{array} {c} \frac{jn}{2} \ \ \  \text{for} \  j \ \text{even,} \ \ j\ne 2  \\
\
\\ i + \frac{(j-3)n}{2} \  \text{for} \  j \ \text{odd,} \ \ j\le \frac{k+1}{2}
 \\ \
 \\ i + \frac{(k+j+\varepsilon)n}{2} \  \text{for} \  j \ \text{odd,} \ j > \frac{k+1}{2} \end{array}  \right.$$

where $\epsilon = -2$ or $0$ according to whether $k\equiv 1$ or $3$ (mod$4$), respectively.

The set ${\cal F} = \{B_1,\dots, B_n\}$ is a $(K_{2nk+1},C_k)$-difference system.

If $n\equiv 2,3$ (mod $4$), we fix a hooked Skolem sequence $(s_1, \dots, s_n)$ of order $n$. In this case, a  difference system is obtained by replacing the cycle  $B_1$ of ${\cal F}$ with the either the cycle $A$, or the cycle $A'$,  according to whether $k\equiv 1$ (mod $4$) or $3 < k \equiv 3$ (mod $4$). Here  $A$ is defined as follows: 
\[
\begin{aligned}
&a_1 =1 \quad a_3 = 0 \quad a_{k-2}=5n \quad a_{k-1}=(4-k)n \quad a_k=(k+3)n/2 +1\\
&a_j=b_{1j} \ \ \ \text{for $j\notin \{1,3,k-2,k-1,k\}$.}
\end{aligned}
\]
%$a_1=1$, $a_3=0$, $a_{k-2}=5n$, $a_{k-1}=(4-k)n$, $a_k=(k+3)n/2 +1$, $a_j=b_{1j}$, for all $j$'s $\notin \{1,3,k-2,k-1,k\}$. W
While $A'$ is defined as follows:
\[
\begin{aligned}
&a_1 =1 \quad a_3 = 0 \\
&a_j=b_{1j} \ \ \ \text{for  $1\ne j \ne 3$.}
\end{aligned}
\] % $a_1=1$, $a_3=0$, $a_j=b_{1j}$ for $1\ne j \ne 3$.
 
We can prove the following

\begin{proposition}\label{Genk}
Let $NC(2nk+1,k)$ be the number of cyclic non-isomorphic $(K_{2nk+1},C_k)$-designs, $k>5$ odd, $n >1$. Therefore:

\begin{itemize}

\item If $n\not\equiv 4$ (mod $7$) is sufficiently large, then
$NC(2nk+1,k)\ge \dfrac{(2.49)^n}{k}$.

\item If $n \equiv 4$ (mod $7$) is sufficiently large, then
$NC(2nk+1,k)\ge \dfrac{(3.35)^n}{k}$.

\end{itemize}
\end{proposition}

\noindent
{\em Proof.} Let $v=2nk+1$ and let ${\cal F }$ be the $(K_v,C_k)$-difference system obtained in \cite {BD1} and above described.  Obviously, all the cycles of ${\cal F}$ are of type $1$ since $1$ is the unique common divisor of $2nk+1$ and $k$. 

Observe also that if $S$ and $S'$ are two different Skolem sequences, respectively hooked Skolem sequences, of the same  order $n$, and if ${\cal F}$ and ${\cal F}$' are the related difference systems giving rise to the cyclic $(K_v,C_k)$-designs,
$D$ and $D'$ respectively, then $D\ne D'$ since $\Lambda ({\cal F}) \ne \Lambda ({\cal F}')$.  Recalling Proposition \ref{L2}, we conclude that $NC(2nk+1,k)\ge 2^n \times Sk(n)$, where $Sk(n)$ is the number of distinct Skolem sequences, respectively hooked Skolem sequences, of order $n$ according to whether $n\equiv 0,1$ (mod $4$) or $n\equiv 2,3$ (mod $4$).  We known that $Sk(n)\ge 2^{\lfloor n/3 \rfloor}$ whenever $n\not \equiv 4$ (mod $7$), while $Sk(n)\ge (6.492)^{(2n-1)/7}$ whenever $n\equiv 4$ (mod $7$). In the first case, we obtain: $NC(2nk+1,k)\ge \dfrac{2^n\cdot 2^{\lfloor n/3 \rfloor -1}}{kn}$ and a numeric calculation shows that this number is bigger than $\dfrac{(2.49)^n}{k}$ for $n\ge 640$.

In the second case, we obtain $NC(2nk+1,k)\ge \dfrac{2^n \cdot (6.492)^{(2n-1)/7}}{2nk}$ and a numeric calculation shows that this number is bigger than $\dfrac{(3.35)^n}{k}$ for $n\ge 375$.  $\qed$

\vskip0.3truecm\noindent
Now we consider the case of a cyclic $(K_{mk},C_k)$-design, with both $m$ and $k$ odd, with $m=2n+1$, $k=2h+1$, $n > 1$, $k >5$. 

We can state the following 

\begin{proposition}\label{Genk2}
Let $NC(2nk+k,k)$ be the number of non-isomorphic cyclic $(K_{2nk+k},C_k)$-designs, with $n >1$ and $k>5$ odd. Therefore:

\begin{itemize}

\item If $n\equiv 3, 5, 24, 26 \  (mod \ 28)$ is sufficiently large,
\par\noindent
then $NC(2nk+k,k)\ge \frac{(3.35)^n}{k}$.

\item If $n\not\equiv 3, 5, 24, 26 \  (mod \ 28)$ is sufficiently large,
\par\noindent
then $NC(2nk+k,k)\ge \frac{(2.49)^n}{k}$.

\end{itemize}
\end{proposition}
\noindent
{\em Proof.} If $k\ne 15$ and $k\ne q$, with $q$ a non prime and a prime power, the way to construct a starter family ${\cal F'}$ is given in \cite{BD1}, as already described for the case $k=5$ in the previous section 3.2. 

More precisely,  let $H=m\Bbb Z_{mk}=\{0,m,2m,\dots, (k-1)m\}$ be the subgroup of the additive group $\Bbb Z_{mk}$. Consider the multipartite graph $K_{m\times k}$ whose parts are the cosets $H+t$, $0\le t\le m-1$. A difference system ${\cal F'}= \{B_1, \dots, B_n\}$ whose orbit under $\Bbb Z_{mk}$ gives a cyclic $(K_{m\times k},C_k)$-design is constructed by choosing each $B_i$ with the following rules:

$$b_{ij}: \left\{\begin{array} {c} \dfrac{m(j-1)}{2} \ \ \  \text{for} \  j \ \text{odd,} \ \ j\ne k  \\ m(h - \dfrac{j}{2})-i \ \ \  \text{for} \  j \ \text{even,} \end{array}  \right.$$

$$b_{ik}=r_i +(hm-n-1),$$

\noindent
where $(r_1,\dots, r_n)$ is a either a split Skolem sequence or a hooked split Skolem sequence of order $n$ according to whether $n\equiv 0,3$ or $n\equiv 1,2$ (mod $4$).

The set ${\cal F'}$ is a difference system for a cyclic $(K_{k\times m},C_k)$ design as $\partial {\cal F'} = \Bbb Z_{km} - m\Bbb Z_{km}$, Theorem 3.2 of \cite{BD1}. Moreover,  all the cycles of ${\cal F'}$ are of type 1. In fact, since $\partial {\cal F'} = \Bbb Z_{km} - m\Bbb Z_{km}$, each cycle of ${\cal F'}$ gives rise to $2k$ distinct differences of $\Bbb Z_{km}$.

If $k\ne 15$ and $k\ne q$, with $q$ a non prime and a prime power, a cyclic $(K_k,C_k)$-design exists (see \cite{BD1} if $k$ is a prime, and \cite{BD2} if $k$ is not a prime) and denote by ${\cal B}$ a  $(K_k,C_k)$-difference system. Applying Theorem 4.1 of \cite{BD1} we know that ${\cal F}= {\cal F'}\cup {\cal B}$ is a $(K_{mk},C_k)$-difference system. Moreover, since the cyclic $(K_k,C_k)$-design is Hamiltonian, each cycle of ${\cal B}$ is of type greater than $1$. If we change the split Skolem sequence which defines ${\cal F'}$ we obtain a new family ${\cal F"}$ and $\bar{{\cal F}}= {\cal F"}\cup {\cal B}$ is still a $(K_{mk},C_k)$-difference system giving rise to a different design, since $\Delta {\cal F'}\ne \Delta {\cal F"}$.

If $k=q$, with $q$ a non prime and a prime power, a $(K_{mk},C_k)$-difference system ${\cal F}_q$ was constructed in \cite{V}. 
It is ${\cal F}_q={\cal F}'_q\cup {\cal B'}$, where all the cycles of ${\cal B'}$ are of type greater than 1 (see the description, p.303 of \cite{V}), while the cycles of ${\cal F}'_q$ are of type $1$ and are obtained by modifying one cycle of the family ${\cal F'}$ constructed in \cite{BD2} and decribed above, (see p.305 of \cite{V}). Namely, ${\cal F}'_q = ({\cal F'} - \{B_1\})\cup \{B_1'\}$. Here $B_1'=(b'_{11},\dots, b'_{1k})$ is obtained by replacing certain vertices of $B_1=(b_{11},\dots, b_{1k})$, but $b'_{1k}$ which remains equal to $b_{1k}$. This assures that if we change the split Skolem sequence which defines ${\cal F'}$, we obtain a new family ${\cal F}''_q$ instead of ${\cal F}'_q$ and $\bar{{\cal F}_q}= {\cal F}''_q\cup {\cal B'}$ is still a $(K_{mk},C_k)$-difference system giving rise to a different design, since $\Delta{\cal F}'_q\ne \Delta{\cal F}''_q$.

If $k=15$, a $(K_{15m},C_{15})$-difference system ${\cal F}={\cal F}' \cup {\cal B}$ was constructed in Theorem 2.1 of \cite{V}. Here the cycles of ${\cal B}$ are of type greater then $1$, while the cycles of ${\cal F}'$ are of type $1$ and are constructed via a split Skolem sequence $(r_1,\dots, r_n)$ of order $n$.  Namely ${\cal F}'=\{C^*_1,C_2, \dots C_n\}$ with $C_i=(0,6m-i,m,5m-i,2m,4m-i, \dots, 6m,-i,r_i + \frac{13m-1}{2})$ for each $1\le i\le n$ and 
$C^*_1$ is obtained fron $C_1$ by repleacing the vertex $6m-1$ with $-5m+1$ and the vertex $6m$ with $6m-1$.  Obviously, also in this case if we change the split Skolem sequence we obtain a new difference system giving rise to a different design.

Now we are able to say that we can achieve a lower bound for $NC(mk,k)$  by using the same process described in the previous Proposition \ref{Genk}. As above we obtain $NC(mk,k) \ge \frac{2^n \times Sk(n)}{mk}$, where: $m=2n+1$, $n >1$, $k >5$, $k$ odd, and $Sk(n)$ is the number of distinct split Skolem sequences or split hooked Skolem sequences of order $n$, according to $n\equiv 0,3$ or $n\equiv 1,2$ (mod $4$). 

As already observed in Proposition \ref{Gen}, if $n$ is sufficiently large, we have either  $Sk(n)\ge (6.492)^{(2n+1)/7}$ or $Sk(n)\ge (6.492)^{(2n-3)/7}$ according to whether $n\equiv 3,24$ or $n\equiv 5,26$ (mod $28$). Numeric calculations show that both the fractions
$\frac{2^n \times (6.492)^{(2n+1)/7}}{(2n+1)k}$ and $\frac{2^n \times (6.492)^{(2n-3)/7}}{(2n+1)k}$ are bigger than $\frac{(3.35)^n}{k}$ starting from $n\ge 403$.
If $n\not\equiv 3, 5, 24, 26 \  (mod \ 28)$ is sufficiently large, we have $Sk(n) \ge 2^{\lfloor n/3\rfloor}$ and numeric calculations show that the fraction $\frac{2^n \times  2^{\lfloor n/3\rfloor}}{(2n+1)k}$ is bigger than $\frac{(2.49)^n}{k}$ starting from $n\ge 597$.
$\qed$

\end{document}